\documentclass[11pt,a4paper]{article}
\usepackage{pifont}
\usepackage{CJK,epsfig,amsfonts,amsgen,amsmath,amstext,amsbsy,amsopn,amsthm
}
\usepackage{mathrsfs}
\usepackage{float}
\usepackage{epsfig,latexsym,amsfonts,amsmath,amssymb}
\usepackage{color}

\usepackage{comment}
\usepackage{geometry}
\usepackage{calc}
\usepackage{tikz}
\usetikzlibrary{decorations.markings}
\tikzstyle{vertex}=[circle, draw, inner sep=0pt, minimum size=6pt]

\newcommand\qbinom[2]{\genfrac{[}{]}{0pt}{}{#1}{#2}}
\newtheorem{proposition}[subsection]{Proposition}
\newtheorem{theorem}[subsection]{Theorem}
\newtheorem{lemma}[subsection]{Lemma}
\newtheorem{conjecture}[subsection]{Conjecture}

\newtheorem{definition}[subsection]{Definition}
\newtheorem{corollary}[subsection]{Corollary}

\newtheorem{remark}[subsection]{Remark}

\newtheorem{case}[subsection]{Case}
\newtheorem{fact}[subsection]{Fact}
\usepackage{indentfirst} 

\begin{document}
\title{On forbidden poset problems in the linear lattice}
\author{\small 
Jimeng Xiao${}^{1,2,3}$ and Casey Tompkins${}^{4,5}$
 \\[2mm]
\small ${}^1$Department of Applied Mathematics, \\
\small Northwestern Polytechnical University, Xi'an, P.R.~China\\
\small ${}^2$Xi'an-Budapest Joint Research Center for Combinatorics,\\ \small Northwestern Polytechnical University, Xi'an, P.R. China\\
\small ${}^3$Alfr\'ed R\'enyi Institute of Mathematics, Hungarian Academy of Sciences\\
\small ${}^4$Department of Mathematics, Karlsruhe Institute of Technology, Germany\\
\small ${}^5$Discrete Mathematics Group, Institute for Basic Science (IBS), Daejeon, Republic of Korea}

\date{}
\maketitle

\baselineskip 20pt

\date{}
\maketitle
 \vspace{4mm}

\begin{abstract}
 
In this note, we determine the maximum size of a $\{\mathrm{V}_{k}, \Lambda_{l}\}$-free family in the lattice of vector subspaces of a finite vector space both in the non-induced case as well as the induced case, for a large range of parameters $k$ and $l$. These results generalize earlier work by Shahriari and Yu.  We also prove a general LYM-type lemma for the linear lattice which resolves a conjecture of Shahriari and Yu.  
\vskip 0.1in \noindent%
 \\[7pt]
\end{abstract}

\newcommand{\lr}[1]{\langle #1\rangle}

\parindent 17pt
\baselineskip 17pt

\section{Introduction}
Given partially ordered sets (posets) $P$ and $Q$, we say that $P$ is a \emph{subposet} of $Q$ if there exists an injection $\phi:P \to Q$ such that $x \le_P y$ implies $\phi(x) \le_Q \phi(y)$. 
If we also have that $\phi(x) \le_Q \phi(y)$ implies $x \le_P y$, then we say $P$ is an \emph{induced subposet} of $Q$.  
Viewing collections of sets as posets under the inclusion relation, we have the following extremal functions, first introduced by Katona and Tarj\'an \cite{kt}. 
For any collection of finite posets $\mathcal{P}$, let $\mathrm{La}(n, \mathcal{P})$ be the maximum size of a family of subsets of $\{1,2,\ldots,n\}$ which does not contain any $P \in \mathcal{P}$ as a subposet, and let $\mathrm{La}^*(n, \mathcal{P})$ be the maximum size of a family of subsets of $\{1,2,\ldots,n\}$ which does not contain any $P \in \mathcal{P}$ as an induced subposet.  
In the case $\mathcal{P} = \{P\}$ for some poset $P$, we instead write simply $\mathrm{La}(n,P)$ and $\mathrm{La}^*(n,P)$. We denote the sum of the $k$ largest binomial coefficients of the form $\binom{n}{i}$ by $\Sigma(n,k)$.

Let $V$ be an $n$-dimensional vector space over a finite field $\mathbb{F}_q$, where $q$ is a prime power. 
The \emph{linear lattice} of dimension $n$ is the poset of subspaces of $V$ under the inclusion relation.
We denote by $\qbinom{V}{k}_q$ the set of all $k$-dimensional subspaces of $V$ (this set is often referred to as a level of the linear lattice). The number of such subspaces is denoted by the $q$-binomial coefficient $\qbinom{n}{k}_q= \prod_{0 \le i < k} \frac{q^{n - i}-1}{q^{k - i}-1}$. 
When $k = 1$, we write $[n]_q = \qbinom{n}{1}_q$. 
Let $[n]_q! = \prod_{i=1}^{n}[i]_q$. 
Then, it is easy to check that
\[
\qbinom{n}{k}_q=\frac{[n]_q!}{[k]_q![n-k]_q!}.
\]

The general study of forbidden poset problems in the linear lattice was initiated by Ghassan and Shahriari \cite{gs}.   For any collection of finite posets $\mathcal{P}$, let $\mathrm{La}_{q}(n, \mathcal{P})$ be the maximum size of a family of subspaces of $V$ (viewed as a poset under inclusion) which does not contain any $P \in \mathcal{P}$ as a subposet, and let $\mathrm{La}_{q}^*(n, \mathcal{P})$ be the maximum size of a family of subspaces of $V$ which does not contain any $P \in \mathcal{P}$ as an induced subposet. 
We write simply $\mathrm{La}_q(n,P)$ and $\mathrm{La}_q^*(n,P)$ if $\mathcal{P} = \{P\}$ for some poset $P$. We denote by $\Sigma_q(n,k)$ the sum of the $k$-largest $q$-binomial coefficients of the form $\qbinom{n}{i}_q$.

Let $\mathrm{V}$ and $\Lambda$ be the posets on three elements $x$, $y$, $z$ defined by the relations $x, y > z$ and $x, y < z$, respectively. 
In $1983$, Katona and Tarj$\mathrm{\acute{a}}$n \cite{kt} proved the following result.

\begin{theorem}[Katona and Tarj$\mathrm{\acute{a}}$n \cite{kt}]\label{thm_kt}
\[
\mathrm{La}(n,\{\mathrm{V},\Lambda\}) = \mathrm{La}^{*}(n,\{\mathrm{V},\Lambda\}) = 2\binom{n-1}{\lfloor\frac{n-1}{2}\rfloor}.
\]
\end{theorem}
\begin{remark}
In Theorem \ref{thm_kt}, an extremal construction is given by the family
\[
\left\{F : 1\notin F, |F| = \lfloor\frac{n-1}{2}\rfloor \right\}\cup \left\{F\cup\left\{1\right\} : 1\notin F, |F| = \lfloor\frac{n-1}{2}\rfloor\right\}.
\]
\end{remark}

Shahriari and Yu~\cite{sy} showed that in the linear lattice we have the following.
\begin{theorem}[Shahriari and Yu~\cite{sy}] \label{t}
\[
\mathrm{La}_{q}(n,\{\mathrm{V},\Lambda\}) = \qbinom{n}{\lfloor\frac{n}{2}\rfloor}_q.
\]
The extremal construction is either $\qbinom{V}{\lfloor\frac{n}{2}\rfloor}_q$ or $\qbinom{V}{\lceil\frac{n}{2}\rceil}_q$, except in the case $n = 3$ and $q=2$, in which we have two other constructions shown in Figure \ref{fig}.
\end{theorem}

We prove the following induced version of Theorem \ref{t}.

\begin{theorem} \label{thm_v}
\[
\mathrm{La}_{q}^{*}(n,\{\mathrm{V},\Lambda\}) = \qbinom{n}{\lfloor\frac{n}{2}\rfloor}_q.
\]
The extremal construction is either $\qbinom{V}{\lfloor\frac{n}{2}\rfloor}_q$ or $\qbinom{V}{\lceil\frac{n}{2}\rceil}_q$, except in the case $n = 3$ and $q=2$, in which we have two other constructions shown in Figure \ref{fig}.
\end{theorem}

Let $\mathrm{V}_{k}$ denote the poset with elements $x_1, x_2,\ldots, x_k, y$ such that $x_1, x_2, \ldots, x_k > y$, and let $\Lambda_{k}$ denote the same poset but with all relations reversed. In the case when $k$ or $l$ is at least 3 only asymptotic results are known for $\mathrm{La}(n,\{ \mathrm{V}_{k}, \Lambda_{l}\})$ (see~\cite{t} and~\cite{d}).

In the linear lattice, on the other hand, one can prove exact results for larger $k$ and $l$ as well.  Shahriari and Yu~\cite{sy} proved the following.

\begin{theorem}[Shahriari and Yu \cite{sy}]\label{oldnoninduced}
Let $n$ be an even integer, and $k, l$ be two integers such that $k, l \le q$. Then
\[
\mathrm{La}_{q}(n, \{\mathrm{V}_{k}, \Lambda_{l}\}) = \qbinom{n}{\frac{n}{2}}_q,
\]
and the only  $\{\mathrm{V}_{k}, \Lambda_{l}\}$-free family of maximum size is $\qbinom{V}{\frac{n}{2}}_q$. 
\end{theorem}

We extend Theorem \ref{oldnoninduced} by weakening the conditions on $k$ and $l$.

\begin{theorem}\label{betternoninduced}
Let $n$ be an even integer, and $k, l$ be two integers such that $k, l \le q^{\frac{n}{2}}$. Then
\[
\mathrm{La}_{q}(n, \{\mathrm{V}_{k}, \Lambda_{l}\}) = \qbinom{n}{\frac{n}{2}}_q,
\]
and the only  $\{\mathrm{V}_{k}, \Lambda_{l}\}$-free family of maximum size is $\qbinom{V}{\frac{n}{2}}_q$. 
\end{theorem}

In the induced case, we have the following two results. 
\begin{theorem}\label{thm_even}
Let $n$ be an even integer and let $k,l$ be two integers such that $k,l \le q$, then
\[
\mathrm{La}^{*}_{q}(n, \{\mathrm{V}_{k}, \Lambda_{l}\}) = \qbinom{n}{\frac{n}{2}}_q,
\]
and the only maximum size  $\{\mathrm{V}_{k}, \Lambda_{l}\}$-free family is $\qbinom{V}{\frac{n}{2}}_q$. 
\end{theorem}

\begin{theorem}\label{thm_odd}
Let $n$ be an odd integer and let $k,l$ be two integers such that $k,l \le (1-\frac{\sqrt{2}}{2})q$, then
\[
\mathrm{La}^{*}_{q}(n, \{\mathrm{V}_{k}, \Lambda_{l}\}) = \qbinom{n}{\frac{n-1}{2}}_q,
\]
and any maximum size  $\{\mathrm{V}_{k}, \Lambda_{l}\}$-free family is either $\qbinom{V}{\frac{n-1}{2}}_q$ or $\qbinom{V}{\frac{n+1}{2}}_q$. 
\end{theorem}

The butterfly poset, $B$, is defined by $4$ elements $a,b,c,d$ with $a,b < c,d$. De Bonis, Katona and Swanepoel \cite{dk} proved the following theorem. 

\begin{theorem}[De Bonis, Katona and Swanepoel \cite{dk}]\label{thm_dk}
\[
\mathrm{La}(n, B) = \Sigma(n,2).
\]
Equality occurs only for a family consisting of the union of two consecutive levels in the Boolean lattice of largest size.
\end{theorem}

We denote by $Y_{k}$ the poset with elements $x_1, x_2, \ldots, x_k, y, z$ such that $x_1 \le x_2 \le \cdots \le x_k \le y, z$ and $Y'_k$ the same poset but with all relations reversed. In the proof of Theorem \ref{thm_dk}, De Bonis, Katona and Swanepoel actually proved a stronger result by determining $\mathrm{La}(n, \{Y_2,Y_2'\})$.  Later pairs of posets $\{Y_k,Y_k'\}$ were investigated for their own sake.  Methuku and Tompkins~\cite{mt} obtained the following theorem.
\begin{theorem}[Methuku and Tompkins \cite{mt}]\label{thm_mt}
Let $k \ge 2$ and $n \ge k+1$, then
\[
\mathrm{La}(n, \{Y_k, Y^{\prime}_k\})=\Sigma(n,k).
\]
\end{theorem}
Martin et al.~\cite{mm} and Tompkins and Wang~\cite{tw}  (these results were strengthened in Gerbner et al.~\cite{gmn}) proved the induced version of Theorem~\ref{thm_mt} independently.

\begin{theorem}[Martin et al.~\cite{mm}, Tompkins and Wang \cite{tw}, Gerbner et al.~\cite{gmn}]\label{thm_mm}
Let $k \ge 2$ and $n \ge k+1$, then
\[
\mathrm{La}^{*}(n, \{Y_k, Y^{\prime}_k\})=\Sigma(n,k).
\]
\end{theorem}

In the vector space setting, Shahriari and Yu~\cite{sy} proved a version of Theorem~\ref{thm_dk} holds.  Namely, they proved the following.
\begin{theorem}[Shahriari and Yu~\cite{sy}]\label{vsp} Let $n \ge 3$ be an integer and $q$ be a power of a prime, then
\[
\mathrm{La}_q(n, B)=\mathrm{La}_q(n, \{Y_2,Y_2'\}) = \Sigma_q(n,2).
\]
Equality occurs only for a family consisting of the union of two consecutive levels in the linear lattice of maximum size. 
\end{theorem}
Furthermore, they posed a conjecture for the case when $\{Y_k,Y_k'\}$ is forbidden.

For any poset $P$, let $|P|$ be the size of $P$ and $h(P)$ be the length of the largest chain in $P$. Burcsi and Nagy~\cite{bn} and Gr\'osz, Methuku and Tompkins~\cite{gm} proved the following theorems for any poset $P$ (another result in this direction was obtain by Chen and Li \cite{cl}).

\begin{theorem}[Burcsi and Nagy \cite{bn}]\label{thm_bn}
For any poset $P$, when $n$ is sufficiently large, we have
\[
\mathrm{La}(n, P) \le \left(\frac{|P|+h(P)}{2}-1\right)\binom{n}{\lfloor\frac{n}{2}\rfloor}.
\]
\end{theorem}

\begin{theorem}[Gr\'osz, Methuku and Tompkins \cite{gm}]\label{thm_gm}
For any poset $P$, when $n$ is sufficiently large, we have
\[
\mathrm{La}(n, P)\le\frac{1}{2^{k-1}}\big(|P|+(3k-5)2^{k-2}(h(P)-1)-1\big)\binom{n}{\lfloor\frac{n}{2}\rfloor},
\]
for any fixed $k$.
\end{theorem}
We will prove that a version of these theorems holds in the vector space case as well.

The rest of this paper is organized as follows. In the next section, we present some preliminary results. Then we will prove Theorems \ref{thm_v}, \ref{thm_even} and \ref{thm_odd} in Section 3. In the last section, we prove a general LYM-type lemma and use this lemma to prove the vector space analogues of Theorems  \ref{thm_mt}, \ref{thm_mm}, \ref{thm_bn} and \ref{thm_gm}. We note that a recent manuscript of Gerbner \cite{gerb} independently initiates a general study of LYM-type properties of the linear lattice and implies some similar results.

\section{Preliminary results}

In this section, let $\mathcal{F}$ be a $\{\mathrm{V}_{k}, \Lambda_{l}\}$-free family of subspaces of $V$, and let $\mathcal{F}_s = \mathcal{F} \cap \qbinom{V}{s}_q$. 
Now, we define the bipartite graph $\left(\mathcal{F}_s \cup \left(\qbinom{V}{s-1}_q \setminus \mathcal{F}_{s-1}\right), E\right)$, where 
\[
E = \left\{\left(A \in \mathcal{F}_s, B \in \left(\qbinom{V}{s-1}_q \setminus \mathcal{F}_{s-1}\right)\right) : B \subset A\right\}.
\]
Let $\mathcal{F}^{\prime}_{s}$ be any subset of $\mathcal{F}_{s}$, and
\[
N_{s-1}(\mathcal{F}^{\prime}_{s}) =\left\{B \in \left(\qbinom{V}{s-1}_q \setminus \mathcal{F}_{s-1}\right) : (A,B)\in E \text{~for some~} A\in \mathcal{F}'_{s}\right\}.
\]

Before beginning the proof, we need some preliminary results.  Lemma \ref{lemma_hall} and Corollary~\ref{corollary_hall} are motivated by an idea from \cite{kt}.

\begin{lemma}\label{lemma_hall}
Let $n$ be an even integer, and let $k,l$ be two integers such that $k,l \le q^{\frac{n}{2}}$.  Then, $\mathrm{La}_{q}(n,\{ \mathrm{V}_{k}, \Lambda_{l}\})$ can be realized with a family $\mathcal{G}$ of subspaces $G$ satisfying $\dim(G) \le \frac{n}{2}$.
\end{lemma}

\noindent\textbf{Proof:~~}
We first prove that for $s \ge \frac{n}{2} + 1$, the bipartite graph $\left(\mathcal{F}_s \cup (\qbinom{V}{s-1} \setminus \mathcal{F}_{s-1}), E\right)$ contains a matching such that every element of $\mathcal{F}_s$ is contained in some edge.  To prove this, it is enough to check the condition of Hall's theorem, that is
\[
|N_{s-1}(\mathcal{F}^{\prime}_{s})|  \ge |\mathcal{F}^{\prime}_{s}|,
\]
for any $\mathcal{F}^{\prime}_{s} \subseteq \mathcal{F}_{s}$.

Since $\mathcal{F}$ is $\Lambda_l$-free, every $s$-dimensional subspace in $\mathcal{F}^{\prime}_s$ has at most $(l-1)$ subspaces in $\mathcal{F}_{s-1}$.  
Hence, every $s$-dimensional subspace in $\mathcal{F}^{\prime}_s$ has at least $[s]_q-l+1$ subspaces in $N_{s-1}(\mathcal{F}^{\prime}_s)$. 
On the other hand, every $(s-1)$-dimensional subspace in $N_{s-1}(\mathcal{F}^{\prime}_s)$ has at most $[n-s+1]_q$ superspaces in $\mathcal{F}^{\prime}_{s}$. We have
\begin{align*}
\frac{|N_{s-1}(\mathcal{F}^{\prime}_s)|}{|\mathcal{F}^{\prime}_s|}
& \geq \frac{[s]_q-l+1}{[n-s+1]_q} 
\geq \frac{\frac{q^s-1}{q-1}-l+1}{\frac{q^{n-s+1}-1}{q-1}}\\
& \ge \frac{q^{\frac{n}{2}+1}-1-(q^{\frac{n}{2}}-1)(q-1)}{q^{\frac{n}{2}}-1} \\
& \geq \frac{q^{\frac{n}{2}}+q-2}{q^{\frac{n}{2}}-1}
\geq 1,
\end{align*}
since $s \geq \frac{n}{2}+1$, $l \leq q^{\frac{n}{2}}$ and $q \ge 2$.
Applying Hall's theorem, let $M$ be a matching which saturates every vertex in $\mathcal{F}_s$, and let $\mathcal{F}^{*}_{s-1}$ be the set of neighbors of $\mathcal{F}_s$ contained in edges of $M$. 
Clearly, $\mathcal{F}^{*}_{s-1} \cap \mathcal{F}_{s-1} = \emptyset$, and $|\mathcal{F}^{*}_{s-1}| = |\mathcal{F}_{s}|$. 

Now, let $t = t(\mathcal{F})$ be the largest integer $s$ satisfying $\mathcal{F}_{s} \not= \emptyset$ in the family $\mathcal{F}$. 
We iteratively replace $\mathcal{F}$ with $(\mathcal{F} \setminus \mathcal{F}_t) \cup \mathcal{F}^{*}_{t-1}$ until $t \le \frac{n}{2}$. Call the resulting family $\mathcal{G}$.
Clearly, $|\mathcal{G}|  = |\mathcal{F}|$, and $\mathcal{G}$ is $\{\mathrm{V}_{k}, \Lambda_{l}\}$-free since $\mathcal{F}$ is  $\{\mathrm{V}_{k}, \Lambda_{l}\}$-free.  $\hfill \Box$

Since linear lattices are symmetric, one can use the same idea to prove the following corollary.

\begin{corollary}\label{corollary_hall}
Let $n$ be an even integer, and $k,l$ be two integers such that $k,l \le q^{\frac{n}{2}}$.  Then, $\mathrm{La}_{q}(n, \{\mathrm{V}_{k}, \Lambda_{l}\})$ can be realized with a family $\mathcal{G}$ of subspaces $G$ satisfying $\dim(G) \ge \frac{n}{2}$.
\end{corollary}

The next technical lemma will be needed for determining the structure of the extremal families.

\begin{lemma}\label{lemma_structure}
Let $\mathcal{F}$ be a family such that $\dim(F) = \lceil \frac{n}{2}\rceil$ or $\lceil \frac{n}{2}\rceil+1$ for every $F \in \mathcal{F}$. If $|\mathcal{F}| = \qbinom{n}{\lceil \frac{n}{2}\rceil}_q$ and $\mathcal{F}_{\lceil \frac{n}{2}\rceil+1}\neq \emptyset$,  then $\mathcal{F}$ contains a copy of $\Lambda_{q^{\lceil \frac{n}{2}\rceil}}$.
\end{lemma}
\noindent{\bf Proof:~} Any $F \in \mathcal{F}_{\lceil \frac{n}{2}\rceil+1}$ has $[\lceil \frac{n}{2}\rceil+1]_q$ subspaces in $\qbinom{V}{\lceil \frac{n}{2}\rceil}_q$, 
and any $F' \in \qbinom{V}{\lceil \frac{n}{2}\rceil}_q$ has $[\lfloor \frac{n}{2}\rfloor]_q$ superspaces in $\qbinom{V}{\lceil \frac{n}{2}\rceil+1}_q$.
We may now show by a simple averaging argument that there exists an $F \in \mathcal{F}_{\lceil \frac{n}{2}\rceil+1}$ such that $F$ has at least $q^{\lceil \frac{n}{2}\rceil}$ subspaces in $\mathcal{F}_{\lceil\frac{n}{2}\rceil}$.  
Indeed, by the assumption that $|\mathcal{F}| = \qbinom{n}{\lceil \frac{n}{2}\rceil}_q$, the number of relations between $\mathcal{F}_{\lceil \frac{n}{2}\rceil+1}$ and $\mathcal{F}_{\lceil \frac{n}{2}\rceil}$ is at least 
\begin{align*}
  &[\lceil \frac{n}{2}\rceil+1]_q 
\left|\mathcal{F}_{\lceil \frac{n}{2}\rceil+1}\right| 
- [\lfloor \frac{n}{2}\rfloor]_q
 \left|\qbinom{V}{\lceil \frac{n}{2}\rceil}_q
 \setminus\mathcal{F}_{\lceil \frac{n}{2}\rceil}\right|\\
 =& \left|\mathcal{F}_{\lceil \frac{n}{2}\rceil+1}\right|([\lceil \frac{n}{2}\rceil+1]_q-[\lfloor \frac{n}{2}\rfloor]_q)\\
 = & \left|\mathcal{F}_{\lceil \frac{n}{2}\rceil+1}\right|\frac{q^{\lceil \frac{n}{2}\rceil+1}-q^{\lfloor \frac{n}{2}\rfloor}}{q-1}\\
\ge& \left|\mathcal{F}_{\lceil \frac{n}{2}\rceil+1}\right|q^{\lceil \frac{n}{2}\rceil}.
\end{align*} 
Thus, on average an element of $\mathcal{F}_{\lceil \frac{n}{2}\rceil+1}$ contains at least $q^{\lceil \frac{n}{2}\rceil}$ subspaces in $\mathcal{F}_{\lceil \frac{n}{2}\rceil}$. $\hfill \Box$

In the same way one can show the following.
\begin{lemma}\label{lemma_structure_reverse}
Let $\mathcal{F}$ be a family such that $\dim(F) = \lfloor \frac{n}{2}\rfloor$ or $\lfloor \frac{n}{2}\rfloor-1$ for every $F \in \mathcal{F}$. If $|\mathcal{F}| = \qbinom{n}{\lfloor \frac{n}{2}\rfloor}_q$ and $\mathcal{F}_{\lfloor \frac{n}{2}\rfloor-1}\neq \emptyset$,  then $\mathcal{F}$ contains a copy of $\mathrm{V}_{q^{\lceil \frac{n}{2}\rceil}}$.
\end{lemma}
Now, we can prove Theorem \ref{betternoninduced}.  

\noindent\textbf{Proof:~~} 
Combining Lemma \ref{lemma_hall} and Corollary \ref{corollary_hall}, it is easy to see that 
\[
\mathrm{La}_{q}(n, \{\mathrm{V}_{k}, \Lambda_{l}\}) = \qbinom{n}{\frac{n}{2}}_q.
\]

Now we prove that if $\mathcal{F}$ is a $\{\mathrm{V}_{k}, \Lambda_{l}\}$-free family of maximum size, then $\mathcal{F} = \qbinom{V}{\frac{n}{2}}_q$. 
Suppose not. If there is a subspace $F \in \mathcal{F}$ of dimension larger than $\frac{n}{2}$, then we may assume, without loss of generality, that $|\mathcal{F}| = \qbinom{n}{\frac{n}{2}}_q$ and for every $F \in \mathcal{F}$, $\frac{n}{2} \le \dim(F) \le \frac{n}{2} + 1$ and $\mathcal{F}_{ \frac{n}{2}+1}\neq \emptyset$. 
 By Lemma \ref{lemma_structure}, $\mathcal{F}$ contains a copy of $\Lambda_l$, a contradiction. The case when $\mathcal{F}$ contains only subspaces of dimension at most $\frac{n}{2}$ is handled similarly by Lemma~\ref{lemma_structure_reverse}.  $\hfill \Box$

\section{Proofs of Theorems \ref{thm_v}, \ref{thm_even} and \ref{thm_odd}}
In this section, let $\mathcal{F}$ be an induced $\{\mathrm{V}_{k}, \Lambda_{l}\}$-free family, and let $\mathcal{F}_{s} = \qbinom{V}{s}_q \cap \mathcal{F}$. 
We call a subspace $F \in \mathcal{F}$ small if for any other $F' \in \mathcal{F}$, $F' \not\subseteq F$.
For every $A \in \mathcal{F}_{s}$, let $F_1, F_2, \ldots, F_{r}$ be $r$ small proper subspaces  of $A$  in $\mathcal{F}$. 
 We note that the $F_i$'s enumerate \textit{all} of the proper small subspaces of $A$.
Clearly, $0 \le r \le l-1$ since $\mathcal{F}$ is induced $\Lambda_{l}$-free ($r = 0$ if $A$ is small). 
Let $f_1 \subseteq F_1$,  $f_2 \subseteq F_2$,\ldots, $f_{r}\subseteq F_{r}$ be $r$ one dimensional subspaces (note that $f_1, f_2, \ldots, f_{r}$ are not necessarily distinct).
Then, we have the following proposition.
\begin{proposition}\label{prop_small}
If $F$ is a subspace of $A$ and $F \in \mathcal{F}$, then $f_i \subseteq F$ for some $i \in [r]$.
\end{proposition}
\noindent{\bf Proof:~~} The subspace $F$ is either small (suppose $F = F_i$ in this case) or contains some small subspace $F_i$. 
In both cases we have $f_i \subseteq F_i\subseteq F$. $\hfill \Box$

Now, we define a family $M(A)$ collecting all $(s-1)$-subspaces of $A$ which do not contain any of the $f_i$.
\[
M(A) = \left\{B: \dim(B)  =s-1, f_1 \not\subseteq B, f_2 \not\subseteq B, \cdots, f_{r} \not\subseteq B \text{~and~} B \subseteq A\right\}.
\]
By Proposition \ref{prop_small}, we have that the following properties of $M(A)$ hold.
\begin{proposition}\label{prop_ma}
\begin{enumerate}
 \item [\rm(i)] For any $B \in M(A)$ and $F \in \mathcal{F}$ such that $\dim(F) \le s-2$, $F \not\subseteq B$.
\item [\rm(ii)] $M(A) \cap \mathcal{F}_{s-1} = \emptyset.$
\item [\rm(iii)] $|M(A)| \ge [s]_q - (l-1)[s-1]_q.$
\end{enumerate}
\end{proposition}
\noindent\textbf{Proof:~~} \rm(i)  If $F \not\subseteq A$, then $F \not\subseteq B$ since $B \subseteq A$. Let $F \subseteq A$, then by Proposition \ref{prop_small}, $f_i \subseteq F$ for some $i \in [r]$. However, by the definition of $M(A)$, $f_i \not\subseteq B$, and so $F \not\subseteq B.$ 

\rm(ii) Suppose not. Let $B \in M(A) \cap \mathcal{F}_{s-1}$. We have that $B$ contains a one dimensional subspace $f_i$ by Proposition \ref{prop_small}, but $f_i \not\subseteq B$ by the definition of $M(A)$, a contradiction.

\rm(iii) For an $s$-dimensional subspace $A$, there are $[s]_q$ $(s-1)$-dimensional subspaces of $A$. At most $[s-1]_q$ among them contain $f_i$ for each $f_i$. So
\[
|M(A)| \ge  [s]_q - r[s-1]_q \ge [s]_q - (l-1)[s-1]_q,
\] 
as required. $\hfill \Box$

Now, we define the bipartite graph $\left(\mathcal{F}_s \cup \left(\qbinom{V}{s-1}_q \setminus \mathcal{F}_{s-1}\right), E\right)$, where 
\[
E = \left\{(A,B) : A \in \mathcal{F}_s, B \in M(A) \right\}.
\]
Let $\mathcal{F}^{\prime}_{s}$ be any subset of $\mathcal{F}_{s}$, and
\[
N_{s-1}(\mathcal{F}^{\prime}_{s}) =\{B : (A,B)\in E \text{~for some~} A\in \mathcal{F}'_{s}\}.
\]

\begin{lemma}
Let $n$ be an odd integer, and $k,l$ be two integers such that $k, l \le q$. Then, $\mathrm{La}^{*}_{q}(n, \{\mathrm{V}_{k}, \Lambda_{l}\})$ can be realized with a family $\mathcal{G}$ of subspaces $G$ satisfying $\dim(G) \le \frac{n+1}{2}$.
\end{lemma}

\noindent\textbf{Proof:~~} 
We first show that for $s \ge \frac{n+3}{2}$, the bipartite graph $\left(\mathcal{F}_s \cup \left(\qbinom{V}{s-1}_q \setminus \mathcal{F}_{s-1}\right), E\right)$ contains a matching such that every element of $\mathcal{F}_s$ is contained in some edge.  By Hall's theorem, it is enough to prove
\[
|N_{s-1}(\mathcal{F}^{\prime}_{s})|  \ge |\mathcal{F}^{\prime}_{s}|,
\]
for any $\mathcal{F}'_{s} \subseteq \mathcal{F}_s$.

On the one hand, by (ii) from Proposition \ref{prop_ma}, $M(A) \cap \mathcal{F}_{s-1} = \emptyset$, then every $s$-dimensional subspace $A$ in $\mathcal{F}^{\prime}_s$ has $|M(A)|$ subspaces in $N_{s-1}(\mathcal{F}'_{s})$.  
On the other hand, every $(s-1)$-dimensional subspace in $N_{s-1}(\mathcal{F}^{\prime}_s)$ has at most $[n-s+1]_q$ superspaces in $\mathcal{F}^{\prime}_{s}$. Then, by (iii) from Proposition~\ref{prop_ma}, we have
\begin{align*}
\frac{|N_{s-1}(\mathcal{F}^{\prime}_s)|}{|\mathcal{F}^{\prime}_s|}
&\geq \frac{[s]_q-(l-1)[s-1]_q}{[n-s+1]_q}
\geq \frac{\frac{q^s-1}{q-1}-\frac{(l-1)(q^{s-1}-1)}{q-1}}{\frac{q^{n-s+1}-1}{q-1}} \\
&\geq \frac{q^{s-1} + q-2}{q^{n-s+1}-1} 
\geq \frac{q^{\frac{n+1}{2}}+q-2}{q^{\frac{n-1}{2}}-1}\ge 1,
\end{align*}
since $s \geq \frac{n+3}{2}$, $l \leq q$ and $q\ge 2$.
Let $M$ be a matching which saturates every vertex in $\mathcal{F}_s$, and $\mathcal{F}^{*}_{s-1}$ be matched under $M$. 
Clearly, $|\mathcal{F}^{*}_{s-1}| = |\mathcal{F}_{s}|$, and $\mathcal{F}^{*}_{s-1} \cap \mathcal{F}_{s-1} = \emptyset$  by (ii) from Proposition \ref{prop_ma}. 
Now, let $t = t(\mathcal{F})$ be the largest integer $s$ satisfying $\mathcal{F}_{s} \not= \emptyset$ in our family $\mathcal{F}$. 
We repeatedly replace $\mathcal{F}$ by $(\mathcal{F} \setminus \mathcal{F}_t) \cup \mathcal{F}^{*}_{t-1}$ until $t \le \frac{n+1}{2}$.
Call the resulting family $\mathcal{G}$. Clearly, $|\mathcal{G}|  = |\mathcal{F}|$.  
Then it is enough to show that $\mathcal{F}$ is induced $\{\mathrm{V}_{k}, \Lambda_{l}\}$-free in every step.

By contradiction, assume that at some step $\mathcal{F}$ is $\{\mathrm{V}_{k}, \Lambda_{l}\}$-free but  $(\mathcal{F} \setminus \mathcal{F}_t) \cup \mathcal{F}^{*}_{t-1}$ contains $\mathrm{V}_{k}$ or $\Lambda_{l}$.  We distinguish two cases.

\begin{case}
 $(\mathcal{F} \setminus \mathcal{F}_t) \cup \mathcal{F}^{*}_{t-1}$ contains an induced $\Lambda_{l}$.
\end{case}

Let $F_1, F_2, F_3, \ldots, F_l \subset F$ be $l+1$ subspaces in $(\mathcal{F} \setminus \mathcal{F}_t) \cup \mathcal{F}^{*}_{t-1}$ which form  an induced $\Lambda_{l}$. 
Then, $F \in \mathcal{F}^{*}_{t-1}$, since $\mathcal{F}$ is induced $\Lambda_{l}$-free. Let $A$ be matched with $F$ under $M$, then $A$ together with $F_1, F_2, F_3, \ldots, F_l$ form an induced $\Lambda_{l}$ in $\mathcal{F}$, a contradiction. 

\begin{case}
$(\mathcal{F} \setminus \mathcal{F}_t) \cup \mathcal{F}^{*}_{t-1}$ contains an induced $\mathrm{V}_{k}$.
\end{case}

Let $F \subset F_1, F_2, F_3, \ldots, F_k$ be $k+1$ subspaces in $(\mathcal{F} \setminus \mathcal{F}_t) \cup \mathcal{F}^{*}_{t-1}$ which form an induced $\mathrm{V}_{k}$. Since $\mathcal{F}$ is induced $\mathrm{V}_{k}$-free, we may suppose  $F_1 \in \mathcal{F}^{*}_{t-1}.$
Let $A$ be matched with $F_1$ under $M$. Then $F_1 \in M(A)$. Note that $F \in \mathcal{F}$ and $\dim(F) \le t-2$. By (i) from Proposition \ref{prop_ma}, $F \not\subseteq F_1$, a contradiction.  $\hfill \Box$

Using the same ideas from Section 2, one can prove the following corollaries similarly.

\begin{corollary}\label{lemma_induced_hall}
Let $n$ be an odd integer, and $k,l$ be two integers such that $k, l \le q$. Then, $\mathrm{La}^{*}_{q}(n, \{\mathrm{V}_{k}, \Lambda_{l}\})$ can be realized with a family $\mathcal{G}$ of subspaces $G$ satisfying $\frac{n-1}{2}\le\dim(G) \le \frac{n+1}{2}$.
\end{corollary}

\begin{corollary}\label{lemma_induced_hall_even}
Let $n$ be an even integer, and $k,l$ be two integers such that $k, l \le q$. Then, $\mathrm{La}^{*}_{q}(n, \{\mathrm{V}_{k}, \Lambda_{l}\})$ can be realized with a family $\mathcal{G}$ of subspaces $G$ satisfying $\dim(G) = \frac{n}{2}$.
\end{corollary}

Theorem \ref{thm_even} follows from Corollary \ref{lemma_induced_hall_even} and the equality cases are again settled by applying Lemmas \ref{lemma_structure} and \ref{lemma_structure_reverse}. (Once a family is contained in two levels there is no distinction between an induced and noninduced copy of $\mathrm{V}_{k}$ or $\Lambda_{l}$.)

Now, we turn to prove Theorem \ref{thm_odd}. Before beginning the proof, we need the following lemma.

\begin{lemma}\label{lemma_proplane}
Let $V_3$ be a 3-dimensional vector space over $\mathbb{F}_{q}$. 
If $\mathcal{F} \subseteq (\qbinom{V_3}{1}_q\cup\qbinom{V_3}{2}_q)$ is $\{V_{k}, \Lambda_{l}\}$-free, where $k,l \le q-\frac{\sqrt{2}}{2}q$, then
\[
|\mathcal{F}| \leq q^2 + q + 1,
\] 
and the only families which attain equality are $\qbinom{V_3}{1}_q$ and $\qbinom{V_3}{2}_q$.
\end{lemma}

\noindent \textbf{Proof:~~}
Let $\mathcal{F} = \mathcal{A} \cup \mathcal{B}$, where $\mathcal{A} \subseteq \qbinom{V_3}{2}_q$ and
$\mathcal{B} \subseteq \qbinom{V_3}{1}_q$, and let  $\mathcal{A^{\prime}} = \qbinom{V_3}{2}_q \setminus \mathcal{A}$ and $\mathcal{B^{\prime}} = \qbinom{V_3}{1}_q \setminus \mathcal{B}$. 

We prove the inequality by contradiction.
Suppose that $|\mathcal{F}|=|\mathcal{A}|+|\mathcal{B}| \geq q^2+q+2$. Note that $\left|\qbinom{V_3}{1}_q\right| = \left|\qbinom{V_3}{2}_q\right| = q^2 +q+1$, so we have $|\mathcal{A}| > |\mathcal{B^{\prime}}|$ and $|\mathcal{B}| > |\mathcal{A^{\prime}}|$. 
Since $\mathcal{F}$ is $\Lambda_{l}$-free, for every $A \in \mathcal{A}$, the number of subspaces of $A$ in $\mathcal{B}$ is at most $l-1$, thus the number of subspaces of $A$ in $\mathcal{B^{\prime}}$ is at least $(q+1)-(l-1) = q + 2 - l$. 
Since $|\mathcal{A}| > |\mathcal{B^{\prime}}|$, there exists a subspace $B \in \mathcal{B^{\prime}}$ with at least $q+3-l$ superspaces $A_1, A_2, \ldots, A_{q+3-l}$ in $\mathcal{A}$ by the pigeonhole principle.

\begin{figure}[htbp]
\begin{tikzpicture}

\draw [fill] (-5.5,-3.5) circle [radius = 0.05];
\node[left] at(-5.4,-0.75){$A_1$};
\node[left] at(-5.5,-3.3){$B$};
\draw [fill] (-5.5,-0.5) circle [radius = 0.05];
\node[left] at(-4.5,-0.75){$A_2$};
\draw [fill] (-4.5,-0.5) circle [radius = 0.05];
\node[right] at(-0.4,-0.8){$A_{q+3-l}$};
\node[left] at(-2.5,-0.8){$\cdots$};
\node[left] at(-1,-3.2){$\cdots$};
\draw [fill] (-0.5,-0.5) circle [radius = 0.05];
\draw [fill] (-5.0,-3.5) circle [radius = 0.05];
\draw [fill] (-4.8,-3.5) circle [radius = 0.05];
\node[left] at(-4.2,-3.2){$\cdots$};
\node[rotate = 0] at (-4.6, -3.7) {$\underbrace{\hspace{0.8cm}}$};
\node[below] at(-4.6,-3.7){$q+1-l$};
\draw [fill] (-4.2,-3.5) circle [radius = 0.05];
\draw [fill] (-3.1,-3.5) circle [radius = 0.05];
\draw [fill] (-3.3,-3.5) circle [radius = 0.05];
\node[left] at(-2.55,-3.2){$\cdots$};
\node[rotate = 0] at (-2.9, -3.7) {$\underbrace{\hspace{0.8cm}}$};
\node[below] at(-2.9,-3.7){$q+1-l$};
\draw [fill] (-2.5,-3.5) circle [radius = 0.05];
\draw [fill] (0,-3.5) circle [radius = 0.05];
\draw [fill] (0.2,-3.5) circle [radius = 0.05];
\node[left] at(0.8,-3.2){$\cdots$};
\node[rotate = 0] at (0.4, -3.7) {$\underbrace{\hspace{0.8cm}}$};
\node[below] at(0.4,-3.7){$q+1-l$};
\draw [fill] (0.8,-3.5) circle [radius = 0.05];
\node[above] at(-3,-0.3){$q+3-l$};
\node[rotate = 0] at (-3, -0.3) {$\overbrace{\hspace{5cm}}$};
\node[below] at(-7,-0.2){$\qbinom{V_3}{2}_q$};
\node[above] at(-2,0.6){$\mathcal{A}$};
\node[above] at(4,0.6){$\mathcal{A}'$};
\node[rotate = 0] at (-2, 0.5) {$\overbrace{\hspace{8cm}}$};
\node[rotate = 0] at (4, 0.5) {$\overbrace{\hspace{4cm}}$};
\node[below] at(-7,-3.2){$\qbinom{V_3}{1}_q$};
\node[below] at(-2.5,-4.6){$\mathcal{B'}$};
\node[below] at(3.5,-4.6){$\mathcal{B}$};
\node[rotate = 0] at (-2.5, -4.5) {$\underbrace{\hspace{7cm}}$};
\node[rotate = 0] at (3.5, -4.5) {$\underbrace{\hspace{5cm}}$};
\draw[-][thick] (-6,-0.5) -- (6,-0.5); 
\draw[-][thick] (-6,-3.5) -- (6,-3.5); 
\draw[-][thick] (-5.5,-0.5) -- (-5.5,-3.5); 
\draw[-][thick] (-5.5,-0.5) -- (-5.0,-3.5); 
\draw[-][thick] (-5.5,-0.5) -- (-4.8,-3.5); 
\draw[-][thick] (-5.5,-0.5) -- (-4.2,-3.5); 
\draw[-][thick] (-4.5,-0.5) -- (-5.5,-3.5); 
\draw[-][thick] (-4.5,-0.5) -- (-3.3,-3.5); 
\draw[-][thick] (-4.5,-0.5) -- (-3.1,-3.5); 
\draw[-][thick] (-4.5,-0.5) -- (-2.5,-3.5); 
\draw[-][thick] (-0.5,-0.5) -- (-5.5,-3.5); 
\draw[-][thick] (-0.5,-0.5) -- (0,-3.5); 
\draw[-][thick] (-0.5,-0.5) -- (0.2,-3.5); 
\draw[-][thick] (-0.5,-0.5) -- (0.8,-3.5); 
\end{tikzpicture}
\caption{$\qbinom{V_3}{1}_q\cup\qbinom{V_3}{2}_q$.}
\end{figure}

For $1 \le i,j \le q+3-l$, $A_i$ and $A_j$ have only one common subspace $B$, since there is no butterfly in two consecutive levels of a linear lattice. So we have 
\[
|\mathcal{B^{\prime}}| 
\geq (q+3-l)(q+1-l)+1,
\]
and similarly, we have
\[
|\mathcal{B}| 
> |\mathcal{A^{\prime}}| 
\geq (q+3-k)(q+1-k)+1,
\]
since $\mathcal{F}$ is $\mathrm{V}_{k}$-free. 
Then,
\[
q^{2} + q +1 = \left|\qbinom{V}{1}_q\right|=|\mathcal{B^{\prime}}| + |\mathcal{B}| 
> (q+3-l)(q+1-l) + (q+3-k)(q+1-k)+2,
\]
a contradiction when $k,l \le q-\frac{\sqrt{2}}{2}q$. This completes the proof of the inequality.
Furthermore, if $|\mathcal{F}| = q^2+q+1$ and $\mathcal{A}, \mathcal{B} \not= \emptyset$, we will have $|\mathcal{A}| = |\mathcal{B'}|$ instead of $|\mathcal{A}| > |\mathcal{B'}|$. Then there exists a subspace $B \in \mathcal{B^{\prime}}$ with at least $q+2-l$ superspaces in $\mathcal{A}$, and so
\[
q^{2} + q +1 
= |\mathcal{F}| \ge (q+2-l)(q+1-l) + (q+2-k)(q+1-k)+2,
\]
but this contradicts the condition $k,l \le q-\frac{\sqrt{2}}{2}q$. $\hfill \Box$

\begin{remark}\label{rem}
In Lemma \ref{lemma_proplane}, the upper bound of $|\mathcal{F}|$ is true when the weaker condition $q^{2} + q +1 < (q+3-l)(q+1-l) + (q+3-k)(q+1-k)+2$ holds, and the extremal structure of $\mathcal{F}$ holds when the weaker condition $q^{2} + q +1 < (q+2-l)(q+1-l) + (q+2-k)(q+1-k)+2$ is satisfied.
\end{remark}

Now, we are ready to prove Theorem \ref{thm_odd}.

\noindent\textbf{Proof of Theorem \ref{thm_odd}:~~}  A maximal chain in a linear lattice of dimension $n$ is a sequence of subspaces $V_0,V_1,\dots,V_n$ where $\{0\}=V_0 \subset V_1 \subset \dots \subset V_n=V$.  We denote by $\mathcal{C}$ the set of all maximal chains in a linear lattice. Now, we double count the number of pairs $(F,C)$, where $F\in \mathcal{F}$, $C\in \mathcal{C}$ such that $F$ is in the chain $C$. 

For every $F\in \mathcal{F}$, there are $[\dim(F)]_q![n-\dim(F)]_q!$ maximal chains though $F$. 
On the other hand,  by Corollary \ref{lemma_induced_hall}, we may assume that for every $F \in \mathcal{F}$, $(n-1)/2 \le \dim(F) \le (n+1)/2$. Then,
we consider a pair of subspaces $(G_{1}, G_{2})$ such that $\dim(G_{1})=\frac{n+3}{2}$, $\dim(G_{2})=\frac{n-3}{2}$ and $G_{2} \subseteq G_{1}$. The subfamily of  $\mathcal{F}$ between $G_1$ and $G_2$ satisfies the condition of Lemma \ref{lemma_proplane}. Hence, the size of the subfamily can be bounded as $q^2 + q+ 1$,
and the number of chains between $G_1$ and $G_2$ though some $F$  in the subfamily is $(q^2 + q+ 1)(q+1)$.
Clearly, the number of maximal chains between $\{0\}$ and $G_2$ ($G_1$ and $V$) is $[\frac{n-3}{2}]_q!$, and the number of such pairs $(G_{1}, G_{2})$ is $\qbinom{n}{\frac{n+3}{2}}_q\qbinom{\frac{n+3}{2}}{\frac{n-3}{2}}_q$. 
Then, we have
\begin{equation} \label{ineq}
\sum_{F\in \mathcal{F}}[\dim(F)]_q![n-\dim(F)]_q!
\leq\qbinom{n}{\frac{n+3}{2}}_q\qbinom{\frac{n+3}{2}}{\frac{n-3}{2}}_q([\frac{n-3}{2}]_q!)^2(q^2+q+1)(q+1)
=[n]_q!.
\end{equation}
It follows that
\[
\frac{|\mathcal{F}|}{\qbinom{n}{\frac{n+1}{2}}_q}
\leq\sum_{F\in \mathcal{F}} \frac{1}{\qbinom{n}{\dim(F)}_q}
\leq1,
\]
hence, 
\[
|\mathcal{F}| \leq \qbinom{n}{\frac{n+1}{2}}_q.
\]
This completes the proof of $\mathrm{La}^*_{q}(n, \{\mathrm{V}_{k}, \Lambda_{l}\}) = \qbinom{n}{\frac{n+1}{2}}_q$.  Since equality must hold in the first inequality of~\eqref{ineq} when $|\mathcal{F}| = \qbinom{n}{\frac{n+1}{2}}_q$, we have the following.
\begin{fact}\label{fact_equality}
If $\frac{n-1}{2} \le \dim(F) \le \frac{n+1}{2}$ for every $F\in \mathcal{F}$ and $|\mathcal{F}| = \qbinom{n}{\frac{n+1}{2}}_q$, then the size of the subfamily between $G_1$ and $G_2$ is $q^2+q+1$ for any pair $(G_1,G_2)$ such that $\dim(G_{1})=\frac{n+3}{2}$, $\dim(G_{2})=\frac{n-3}{2}$ and $G_{2} \subseteq G_{1}$.
\end{fact}

We will also make use of the following simple lemma.
\begin{lemma}\label{reviewer}
Let $G$ be a connected, regular, bipartite graph with parts $A$ and $B$.  If for some $A' \subset A$ we have $|N(A')|=|A|$, then either $A'=\emptyset$ or $A'=A$.  
\end{lemma}

Now, we show the largest  induced $\{\mathrm{V}_{k}, \Lambda_{l}\}$-free family is either $\qbinom{V}{\frac{n+1}{2}}_q$ or $\qbinom{V}{\frac{n-1}{2}}_q$ by considering three cases.  Assume $\mathcal{F}$ is an induced $\{\mathrm{V}_{k}, \Lambda_{l}\}$-free family of size  $\qbinom{n}{\frac{n+1}{2}}_q$.

\begin{case}
For all $F\in \mathcal{F}$, $\dim(F) \ge  \frac{n+1}{2}$.
\end{case}
We will show that $\mathcal{F}$ contains only subspaces of dimension $\frac{n+1}{2}$. 
Suppose not, then we may find an induced $\{\mathrm{V}_{k}, \Lambda_{l}\}$-free family $\mathcal{F}'$ of size  $\qbinom{n}{\frac{n+1}{2}}_q$ containing only subspaces of dimension $\frac{n+1}{2}$ and $\frac{n+3}{2}$ such that $\mathcal{F}'_{\frac{n+3}{2}} \neq \emptyset$.
However, by Lemma~\ref{lemma_structure}, $\mathcal{F}$ contains a copy of $\Lambda_l$, since $l \le q - \frac{\sqrt{2}}{2}q \le q^{\frac{n+1}{2}}$, a contradiction. The following case can be proved by a similar argument.

 \begin{case}
For all $F\in \mathcal{F}$, $\dim(F) \le  \frac{n-1}{2}$.
\end{case}

\begin{case}
There exist two subspaces $F_1$ and $F_2$ in $\mathcal{F}$ such that $\dim(F_1) \ge  \frac{n+1}{2}$ and $\dim(F_2) \le  \frac{n-1}{2}$.
\end{case}
We will show this case is impossible.  By Corollary~\ref{lemma_induced_hall}, there is an induced $\{\mathrm{V}_{k}, \Lambda_{l}\}$-free family $\mathcal{F}'$ of the same size as $\mathcal{F}$ such that $\frac{n-1}{2} \le \dim(F') \le \frac{n+1}{2}$ for every $F' \in \mathcal{F}'$.
Clearly, $\left|\mathcal{F}'_{\frac{n+1}{2}}\right| + \left|\mathcal{F}'_{\frac{n-1}{2}}\right| = \qbinom{n}{\frac{n-1}{2}}_q$. 
Then, by the assumption on the dimensions of $F_1$ and $F_2$,  we have
$\mathcal{F}'_{\frac{n+1}{2}} \not= \emptyset$ and $\mathcal{F}'_{\frac{n-1}{2}} \not= \emptyset$. 

Let $N(\mathcal{F}_{\frac{n+1}{2}}')$ be the collection of subspaces of dimension $\frac{n-1}{2}$ contained in subspaces in $\mathcal{F}'$.  Double counting the pairs of subspaces $(A,B)$ where $B \in \mathcal{F}_{\frac{n+1}{2}}'$ and $A$ is a $\frac{n-1}{2}$ dimensional subspace of $B$ we have $ |N(\mathcal{F}_{\frac{n+1}{2}}')| \ge |\mathcal{F}_{\frac{n+1}{2}}'|$. 

Suppose $|N(\mathcal{F}_{\frac{n+1}{2}}')| > \left|\mathcal{F}_{\frac{n+1}{2}}'\right|$, then we must have a pair of subspaces in the family related by containment (by our assumption that $\mathcal{F}'$ has size $\qbinom{n}{\frac{n+1}{2}}_q$),  but this contradicts Fact~\ref{fact_equality} and Lemma~\ref{lemma_proplane}. 

Now suppose $|N(\mathcal{F}_{\frac{n+1}{2}}')| = \left|\mathcal{F}_{\frac{n+1}{2}}'\right|$, then it follows from Lemma~\ref{reviewer} that $\mathcal{F}_{\frac{n+1}{2}}'$ is either empty or the complete level.  This completes the proof of Theorem~\ref{thm_odd}.$\hfill \Box$

\begin{remark}\label{rem_doublecount}
Remark \ref{rem} also holds for Theorem \ref{thm_odd}.
\end{remark}

Now we turn to the proof of Theorem \ref{thm_v}. 
Clearly, the even case follows from Theorem~\ref{thm_even}. So we need to prove the case when $n$ is odd.

\noindent\textbf{Proof of Theorem \ref{thm_v}:~~} 
By Remark \ref{rem_doublecount}, when $k = l = 2$, the weaker condition for upper bound on $|\mathcal{F}|$ is 
$q^2-q-1 >0$.
This is true for $q \ge 2$, and this completes the proof of $\mathrm{La}^*_{q}(n,\{ \mathrm{V}, \Lambda\}) = \qbinom{n}{\frac{n+1}{2}}_q$.
Furthermore, the weaker condition for structure of $\mathcal{F}$ is 
$q^2-3q+1>0$,
and this inequality is true for $q \ge 3$. 
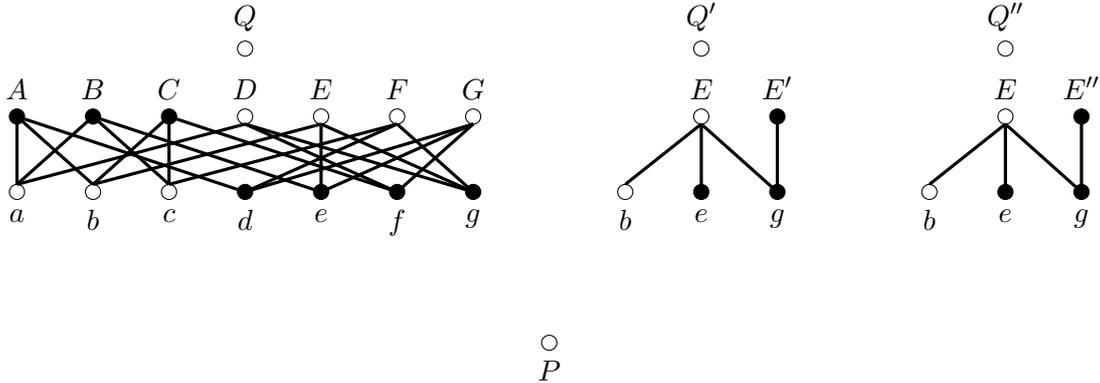
\begin{figure}[htbp]
\begin{tikzpicture}

\draw [fill] (-7,1) circle [radius = 0.1];
\node[above] at(-7,1.1){$A$};

\draw [fill] (-6,1) circle [radius = 0.1];
\node[above] at(-6,1.1){$B$};

\draw [fill] (-5,1) circle [radius = 0.1];
\node[above] at(-5,1.1){$C$};

\draw [] (-4,1) circle [radius = 0.1];
\node[above] at(-4,1.1){$D$};

\draw [] (-3,1) circle [radius = 0.1];
\node[above] at(-3,1.1){$E$};

\draw [] (-2,1) circle [radius = 0.1];
\node[above] at(-2,1.1){$F$};

\draw [] (-1,1) circle [radius = 0.1];
\node[above] at(-1,1.1){$G$};

\draw [] (2,1) circle [radius = 0.1];
\node[above] at(2,1.1){$E$};

\draw [] (6,1) circle [radius = 0.1];
\node[above] at(6,1.1){$E$};

\draw [fill] (3,1) circle [radius = 0.1];
\node[above] at(3,1.1){$E^{\prime}$};

\draw [fill] (7,1) circle [radius = 0.1];
\node[above] at(7,1.1){$E^{\prime\prime}$};

\draw [] (-7,0) circle [radius = 0.1];
\node[below] at(-7,-0.1){$a$};

\draw [] (-6,0) circle [radius = 0.1];
\node[below] at(-6,-0.1){$b$};

\draw [] (-5,0) circle [radius = 0.1];
\node[below] at(-5,-0.1){$c$};

\draw [fill] (-4,0) circle [radius = 0.1];
\node[below] at(-4,-0.1){$d$};

\draw [fill] (-3,0) circle [radius = 0.1];
\node[below] at(-3,-0.1){$e$};

\draw [fill] (-2,0) circle [radius = 0.1];
\node[below] at(-2,-0.1){$f$};

\draw [fill] (-1,0) circle [radius = 0.1];
\node[below] at(-1,-0.1){$g$};

\draw [] (1,0) circle [radius = 0.1];
\node[below] at(1,-0.1){$b$};

\draw [fill] (2,0) circle [radius = 0.1];
\node[below] at(2,-0.1){$e$};

\draw [fill] (3,0) circle [radius = 0.1];
\node[below] at(3,-0.1){$g$};

\draw [] (5,0) circle [radius = 0.1];
\node[below] at(5,-0.1){$b$};

\draw [fill] (6,0) circle [radius = 0.1];
\node[below] at(6,-0.1){$e$};

\draw [fill] (7,0) circle [radius = 0.1];
\node[below] at(7,-0.1){$g$};
\draw [] (-4,1.9) circle [radius = 0.1];
\node[above] at(-4,2){$Q$};
\draw [] (0,-2) circle [radius = 0.1];
\node[below] at(0,-2.1){$P$};

\draw [] (2,1.9) circle [radius = 0.1];
\node[above] at(2,2){$Q^{\prime}$};
\draw [] (6,1.9) circle [radius = 0.1];
\node[above] at(6,2){$Q^{\prime\prime}$};

\draw[-][very thick] (-7,1) -- (-7,0.1); 
\draw[-][very thick] (-7,1) -- (-6,0.1); 
\draw[-][very thick] (-7,1) -- (-4,0); 
\draw[-][very thick] (-6,1) -- (-7,0.1); 
\draw[-][very thick] (-6,1) -- (-5,0.1); 
\draw[-][very thick] (-6,1) -- (-3,0); 
\draw[-][very thick] (-5,1) -- (-6,0.1); 
\draw[-][very thick] (-5,1) -- (-5,0.1); 
\draw[-][very thick] (-5,1) -- (-2,0); 
\draw[-][very thick] (-4,0.9) -- (-7,0.1); 
\draw[-][very thick] (-4,0.9) -- (-2,0); 
\draw[-][very thick] (-4,0.9) -- (-1,0); 
\draw[-][very thick] (-3,0.9) -- (-6,0.1); 
\draw[-][very thick] (-3,0.9) -- (-3,0); 
\draw[-][very thick] (-3,0.9) -- (-1,0); 
\draw[-][very thick] (-2,0.9) -- (-5,0.1); 
\draw[-][very thick] (-2,0.9) -- (-4,0); 
\draw[-][very thick] (-2,0.9) -- (-1,0); 
\draw[-][very thick] (-1,0.9) -- (-4,0); 
\draw[-][very thick] (-1,0.9) -- (-3,0); 
\draw[-][very thick] (-1,0.9) -- (-2,0); 
\draw[-][very thick] (2,0.9) -- (1,0.1); 
\draw[-][very thick] (2,0.9) -- (2,0); 
\draw[-][very thick] (2,0.9) -- (3,0); 
\draw[-][very thick] (6,0.9) -- (5,0.1); 
\draw[-][very thick] (6,0.9) -- (6,0); 
\draw[-][very thick] (6,0.9) -- (7,0); 
\draw[-][very thick] (3,0.9) -- (3,0); 
\draw[-][very thick] (7,0.9) -- (7,0); 
\end{tikzpicture}

\caption{Small examples and illustration of the proof of Theorem \ref{thm_v}.}\label{fig}
\end{figure}

For $q=2$, we can list all the cases for $n= 3$, and there are two constructions which are not levels. (See Figure~\ref{fig}: $\{A,B,C,d,e,f,g\}$ (solid vertices) and $\{ a,b,c,D,E,F,G \}$ (hollow vertices) are the two examples.)
 
Note that in this structure, there is a matching with $3$ edges connecting $6$ subspaces and a single isolated subspace. In Figure~\ref{fig}, $\{A,d\}$, $\{B,e\}$ and $\{C,f\}$ form the matching with $3$ edges,  and $g$ is the single isolated subspace.  (Similarly, $\{D,a\}$, $\{E,b\}$ and $\{F,c\}$ form the matching with $3$ edges, and $G$ is the single isolated subspace.)

However, these constructions do not extend beyond the case $n > 3$ for $q=2$. Similarly, we will prove the only induced $\{\mathrm{V
}, \Lambda\}$-free family of maximum size is either $\qbinom{V}{\frac{n+1}{2}}_q$ or $\qbinom{V}{\frac{n-1}{2}}_q$ by three cases  when  $n > 3$ and $q=2$.
The first two cases $\dim(F) \ge  \frac{n+1}{2}$ for all $F\in \mathcal{F}$ or $\dim(F) \le  \frac{n-1}{2}$ for all $F\in \mathcal{F}$ can be proved by Lemmas  \ref{lemma_structure} and \ref{lemma_structure_reverse}, since $q^{\lceil\frac{n}{2}\rceil} \ge 2$.
Again for the third case, by Corollary \ref{lemma_induced_hall}, we may suppose that for every $F \in \mathcal{F}$, $\frac{n-1}{2} \le \dim(F)
\le \frac{n+1}{2}$, and $\mathcal{F}_\frac{n-1}{2} \not= \emptyset$ and $\mathcal{F}_\frac{n+1}{2} \not=\emptyset$.
Then, again by Lemma~\ref{reviewer} and the assumption that  $|\mathcal{F}|=\qbinom{n}{\frac{n+1}{2}}_q$, we can find two subspaces (say $d \subset A$) in $\mathcal{F}$, where $\dim(d) = \frac{n-1}{2}$ and $\dim(A) = \frac{n+1}{2}$. 

 Instead of using Lemma~\ref{lemma_proplane} to derive a contradiction, we must consider a more detailed argument.
Suppose $\dim(Q) =\frac{n+3}{2}$ and $\dim(P) =\frac{n-3}{2}$ such that $P \subset d \subset A \subset Q$. 
We apply Fact~\ref{fact_equality} for the pair $(Q, P)$. Then we have $7$ (that is, $q^2 + q + 1$) subspaces in $\mathcal{F}$ between $P$ and $Q$.  Since $d$ and $A$ are not in the same level, without loss of generality, we can suppose that  $A,B,C,d,e,f,g \in\mathcal{F}$ and $a,b,c,D,E,F,G \notin \mathcal{F}$ (as in Figure \ref{fig}). 

Since $n \ge 5$, we have $[\frac{n-1}{2}]_q \ge 3$ superspaces of dimension $\frac{n+3 }{2}$ for every subspace of dimension $\frac{n+1}{2}$. 

Thus, $E$ has  two other $\frac{n+3 }{2}$-dimensional superspaces: $Q^{\prime}$ and $Q^{\prime\prime}$.
Now, we apply Fact~\ref{fact_equality} for the pairs $(Q^{\prime}, P)$ and $(Q^{\prime\prime}, P)$.  
Then we have $7$ subspaces in $\mathcal{F}$ between $P$ and $Q^{\prime}$ ($P$ and $Q''$). 
(When we consider a ($\mathrm{V}, \Lambda)$-free family of size $7$ between $P$ and $Q^{\prime}$ ($P$ and $Q''$), the same argument applies as in the $n=3$ case above.) Clearly, there are only $3$ subspaces $b$, $e$ and $g$ of $E$ containing $P$. (See Figure \ref{fig}.) 
So these $3$ subspaces $b$, $e$ and $g$ are also among the $7$ total $(\frac{n-1}{2})$-dimensional subspaces between $P$ and $Q^{\prime}$ ($P$ and $Q''$). 
Since there are $7$ total $(\frac{n-1}{2})$-dimensional subspaces and $7$ total $(\frac{n+1}{2})$-dimensional subspaces between $P$ and $Q^{\prime}$ ($P$ and $Q''$), by the assumption $b \notin\mathcal{F}$ and $e, g \in\mathcal{F}$,  all $7$ subspaces in $\mathcal{F}$ between $P$ and $Q^{\prime}$ ($P$ and $Q''$) cannot form a level, and so they form a matching of $3$ edges and a single isolated subspace.

Now, we will show that $e$ is the single isolated subspace. Suppose not, say $B'$ is the superspace of $e$ in $\mathcal{F}$ between $P$ and $ Q^{\prime}$ ($P$ and $Q''$). 
We have $B \not= B'$, since otherwise $B, E, Q$ and  $Q'$ $(Q'')$ form a butterfly. 
However, by the assumption $B \in\mathcal{F}$, we have that $B'$, $e$ and $B$ will form an induced $\mathrm{V}$, since $\dim(B)=\dim(B')$ and $B \not= B'$. 
Note that $e$ is the single isolated subspace implies that $g$ is in an edge of the matching formed by subspaces in $\mathcal{F}$ between $P$ and $Q^{\prime}$ ($P$ and $Q''$).
Thus, there exist $E^{\prime} \subset Q^{\prime}$ and $E^{\prime\prime} \subset Q^{\prime\prime}$ in $\mathcal{F}$ such that $g \subset E^{\prime}, E^{\prime\prime}$.  
Note that $E' \not= E''$, otherwise $E, E', Q'$ and $Q''$ would form a butterfly.
It follows that $g$, $E'$ and $E''$ form an induced $\mathrm{V}$ in $\mathcal{F}$, a contradiction. $\hfill \Box$

\section{General LYM-type lemma}
Let $V$ be an $n$-dimensional vector space over a finite field $\mathbb{F}_q$, where $q$ is a prime power. Let $\mathcal{H}$ be a family of subspaces of $V$.  We say that $\mathcal{H}$ is simple if there is a basis $\{v_1,v_2,\dots,v_n\}$ of $V$ such that $\mathcal{H} \subseteq \{\text{span}(S):S \in 2^{\{v_1,v_2,\dots,v_n\}}\}$, where span$(S)$ is the subspace spanned by the basis vectors in $S$. For any poset $P$, let $\alpha(\mathcal{H},P)$ denote the maximum size of a $P$-free subfamily of $\mathcal{H}$. We denote by $N_{i}(\mathcal{H})$ the number of $i$-dimensional subspaces in $\mathcal{H}$.  

We now present a general LYM-type lemma.  The proof comes from adopting the methods from \cite{gm} to a vector space setting.

\begin{theorem}\label{lem_lym}
Let $\mathcal{F}$ be a $P$-free family of subspaces of $V$, and let $\mathcal{H}$ be a simple family of subspaces of $V$, then 

\begin{displaymath}
\sum_{F\in \mathcal{F}} \frac{N_{\dim(F)}(\mathcal{H})}{\qbinom{n}{\dim(F)}_q} \le \alpha(\mathcal{H},P).
\end{displaymath}
In particular, if $N_k(\mathcal{H}) = N$ for a given integer $N$ and all $k$, then
\begin{displaymath}
\sum_{F\in \mathcal{F}} \frac{1}{\qbinom{n}{\dim(F)}_q} \le \frac{\alpha(\mathcal{H},P)}{N}.
\end{displaymath}
\end{theorem}

\noindent\textbf{Proof:~~} 
Consider a set $H \in \mathcal{H}$, and without loss of generality assume $H = \text{span}(\{v_1,v_2,\dots,v_r\})$.  Let $F$ be a subspace of $V$ of dimension $r$. Consider maps $\pi$ which replace the basis $\{v_1,v_2,\dots,v_n\}$ of $V$ with an arbitrary basis $\{w_1,w_2,\dots,w_n\}$ of $V$ and assign $\pi(v_i) = w_i$ for $i \in [n]$. For a set $H = \text{span}(\{v_{i_1},v_{i_2},\dots,v_{i_r}\}) \in \mathcal{H}$,  let  $H^\pi = \text{span}(\{\pi(v_{i_1}),\pi(v_{i_2}),\dots,\pi(v_{i_r})\})$, and set $\mathcal{H}^\pi = \{H^\pi:H \in \mathcal{H}\}$.  

We will double count pairs $(F,\pi)$ such that $F \in \mathcal{F}$ and $F \in \mathcal{H}^\pi$.  Suppose $F \in \mathcal{F}$ and $H \in \mathcal{H}^\pi$ both have dimension $r$. The number of $\pi$ such that $F \in \mathcal{H}^\pi$ is

\begin{displaymath}
(q^r-1)(q^r-q)\dots (q^r-q^{r-1}) (q^{n-r}-1)(q^{n-r}-q)\dots (q^{n-r}-q^{n-r-1}).
\end{displaymath}

Observe that if for two distinct $H_1,H_2 \in \mathcal{H}$ we have $F = H_1^{\pi_1}$ and $F = H_2^{\pi_2}$, then $\pi_1 \neq \pi_2$.  It follows that for each $F \in \mathcal{F}$, there are 

\begin{displaymath}
(q^r-1)(q^r-q)\dots (q^r-q^{r-1}) (q^{n-r}-1)(q^{n-r}-q)\dots (q^{n-r}-q^{n-r-1})N_{\dim(F)}(\mathcal{H})
\end{displaymath}
mappings $\pi$ such that $F \in \mathcal{H}^\pi$. Thus, on the one hand, the number of pairs $(F,\pi)$ is 
\begin{displaymath}
\sum_{F \in \mathcal{F}} (q^r-1)(q^r-q)\dots (q^r-q^{r-1}) (q^{n-r}-1)(q^{n-r}-q)\dots (q^{n-r}-q^{n-r-1})N_{\dim(F)}(\mathcal{H}),
\end{displaymath}
or equivalently,
\begin{equation}
\label{generalleft}
\sum_{F \in \mathcal{F}} [\dim(F)]_q![n-\dim(F)]_q!(q-1)^n N_{\dim(F)}(\mathcal{H}).
\end{equation}
Now suppose we fix a mapping $\pi$.  Since $\mathcal{H}$ and $\mathcal{H}^\pi$ are isomorphic as posets with respect to the subspace relation, we have at most $\alpha(\mathcal{H},P)$ many $F \in \mathcal{F}$ such that $F \in \mathcal{H}^\pi$. Since the total number of mappings $\pi$ is 
\begin{displaymath}
(q^n-1)(q^n-q)\dots (q^n-q^{n-1}),
\end{displaymath}
we have an upper bound on the number of pairs $(F,\pi)$ of 
\begin{displaymath}
(q^n-1)(q^n-q)\dots (q^n-q^{n-1}) \alpha(\mathcal{H},P),
\end{displaymath}
or equivalently,
\begin{equation}
\label{generalright}
[n]_q!(q-1)^n \alpha(\mathcal{H},P).
\end{equation}
Combining \eqref{generalleft} and \eqref{generalright}, we have
\begin{displaymath} 
\sum_{F \in \mathcal{F}} [\dim(F)]_q![n-\dim(F)]_q!N_{\dim(F)}(\mathcal{H}) \le [n]_q! \alpha(\mathcal{H},P),
\end{displaymath}
and rearranging yields the desired inequality.

\begin{remark}\label{rem_ind}
The exact same arguments can be carried out to prove the analogous result when we forbid $P$ as an induced subposet.
\end{remark}

Now, we use Lemma \ref{lem_lym} to prove vector space versions of Theorems  \ref{thm_mt}, \ref{thm_mm}, \ref{thm_bn} and \ref{thm_gm}.  We remark that the vector space of Theorem \ref{thm_mt} was conjectured by Shahriari and Yu.

\begin{conjecture}[Shahriari and Yu \cite{sy}]
Let $k \ge 1$ and $n \ge k+1$, then
\[
\mathrm{La}_q(n, \{Y_k, Y^{\prime}_k\})=\Sigma_q(n,k).
\]
\end{conjecture}
We will show that even the induced version of this conjecture holds.
\begin{theorem} \label{pfconj}
Let $k \ge 1$ and $n \ge k+1$, then
\[
\mathrm{La}_q^*(n, \{Y_k, Y^{\prime}_k\})=\Sigma_q(n,k).
\]
\end{theorem}

Let $\{v_1, v_2, \ldots, v_n\}$ be a basis of $V$ and $\mathcal{I}_n$ be the family of subspaces formed by arranging the basis $\{v_1, v_2, \ldots, v_n\}$ in order around a circle and taking those subspaces (excluding $\{0\}$ and $V$) which are spanned by vectors along this cyclic arrangement. 

\begin{lemma}[Tompkins and Wang \cite{tw}]
Let $\alpha^{*}(\mathcal{I}_n,Y_{k}, Y'_{k})$ be the maximum size of an induced $\{Y_k,Y'_k\}$-free subfamily of $\mathcal{I}_n$, then
\[
\alpha^*(\mathcal{I}_n,Y_{k}, Y'_{k}) = kn.
\]
\end{lemma}

\noindent\textbf{Proof of Theorem \ref{pfconj}:} Clearly, $N_{i}(\mathcal{I}_n) = n$ for all $i$. If $\mathcal{F}$ is an induced $\{Y_k,Y'_k\}$-free family such that $\{0\}, V \notin \mathcal{F}$, it follows from Remark \ref{rem_ind} that 
\[
\sum_{F\in \mathcal{F}}\frac{1}{\qbinom{n}{\dim(F)}_q} \le k, 
\]
and so
\[
|\mathcal{F}| \le \Sigma_q(n,k).
\]
Otherwise, if $\{0\}, V \in \mathcal{F}$, we can assume the result is true for $k-1$. Note that the base case $k=1$ is proved by Theorem \ref{thm_v}.
Then it follows that 
\[
|\mathcal{F}|\le \Sigma_q(n,k-1) + 2 \le \Sigma_q(n,k),
\]
since $\mathcal{F}$ is induced $\{Y_{k-1}, Y'_{k-1}\}$-free. 

Now we may assume that $\{0\} \in \mathcal{F}$ and $V \notin \mathcal{F}$. Let $\mathcal{G} = \mathcal{F}\setminus \{\{0\}\}$. If $|\mathcal{G}|\le \Sigma_q(n,k)-1$, then $|\mathcal{F}|\le \Sigma_q(n,k)$. So we may assume $|\mathcal{G}| = \Sigma_q(n,k)$, and $\mathcal{G}$ is a subfamily of the $k$
(or $k + 1$) largest levels in the linear lattice.
If $n \not= k$ modulo $2$, then $\mathcal{G}$ is uniquely determined (i.e., the largest $k$ levels), and it is easy to find an induced $Y_k$ in $\mathcal{F}$. Indeed, we can find an induced $Y_{k-1}$ in  $\mathcal{G}$, which together with $\{0\}$ form an induced $Y_k$. If $n = k$ modulo $2$, then $\mathcal{G}$ is a subfamily of the  $k + 1$ largest levels $\mathcal{L}_{1}, \mathcal{L}_{2}, \ldots, \mathcal{L}_{k+1}$. If $\mathcal{G} \cap \mathcal{L}_{1} = \emptyset$, then we can find an induced $Y_k$ as in the previous case. Otherwise, let $L_1 \in \mathcal{G} \cap \mathcal{L}_{1}$, then find an induced $Y_{k-2}$ in $\mathcal{L}_{2}\cup \cdots \mathcal{L}_{k+1}$ such that $L_1$ is a subspace of every subspace in this $Y_{k-2}$. It is easy to see that $\{0\}$, $L_1$ together with this $Y_{k-2}$ form an induced $Y_k$ in $\mathcal{F}$.$\hfill \Box$

We now recall some other structures which have been used in double counting arguments for forbidden poset problems.

\begin{definition}
Let $\emptyset=A_0 \subset A_1 \subset \dots \subset A_{n-1} \subset A_n = [n]$ be a maximal chain in the $n$-element Boolean lattice.  Then the $k$-\emph{interval chain} defined from this maximal chain is given by $[A_0,A_k] \cup [A_1,A_{k+1}] \cup \dots \cup [A_{n-k},A_n]$.
 A 2-interval chain is called a \emph{double chain}.
\end{definition}
 
\begin{lemma}[Burcsi and Nagy \cite{bn}]
Given a double chain $\mathcal{D}$, we have $N_{i}(\mathcal{D}) = 2$ for every $i$ and 
\[\alpha(\mathcal{D},P) = |P|+h(P)-2.
\]
\end{lemma}
\begin{lemma}[Gr\'osz, Methuku and Tompkins \cite{gm}]\label{intlemma}
Given a $k$-interval chain $\mathcal{C}_{k}$, we have $N_{i}(\mathcal{C}_{k}) = k$ for every $i$ and 
\[\alpha(\mathcal{C}_{k},P) = \frac{k}{2^{k-1}}\big(|P|+(3k-5)2^{k-2}(h(P)-1)-1\big).
\]
\end{lemma}
By Lemma \ref{lem_lym}, the vector space versions of Theorems \ref{thm_bn} and \ref{thm_gm} follows from the above two lemmas, respectively. 
\begin{remark}
In proving the vector space version of Theorem \ref{thm_gm} from the Lemma \ref{intlemma}, we proceed exactly as in the corresponding proof in \cite{gm} replacing each binomial coefficient with the corresponding $q$-binomial, and verify that all the estimates still hold.)
\end{remark}
\begin{theorem}For any poset $P$, when $n$ is sufficiently large, we have
\[
\mathrm{La}_q(n, P) \le \left(\frac{|P|+h(P)}{2}-1\right)\qbinom{n}{\lfloor\frac{n}{2}\rfloor}_q.
\]
\end{theorem}
\begin{theorem}
For any poset $P$, when $n$ is sufficiently large, we have
\[
\mathrm{La}_q(n, P)\le\frac{1}{2^{k-1}}\big(|P|+(3k-5)2^{k-2}(h(P)-1)-1\big)\qbinom{n}{\lfloor\frac{n}{2}\rfloor}_q,
\]
for any fixed $k$.
\end{theorem} 

\section{Acknowledgements}
The authors would like to thank Gyula O.H. Katona for some useful discussions and also the anonymous referee whose remarks, in particular, simplified the proof of Theorem~\ref{thm_odd} (an earlier version of the manuscript invoked a stronger result of Chowdury and Patk\'os~\cite{cp}).
The research was partially supported by the National Research, Development and Innovation Office NKFIH, grant K116769, by CSC (No. 201706290171) and NSFC (No. 11671320).  The research of Tompkins was supported by IBS-R029-C1.

\end{document}